\documentclass[]{elsart}

\usepackage{natbib}
\usepackage{amssymb}

\begin{document}

\begin{frontmatter}

% Title, authors and addresses

% use the thanksref command within \title, \author or \address for footnotes;
% use the corauthref command within \author for corresponding author footnotes;
% use the ead command for the email address,
% and the form \ead[url] for the home page:
% \title{Title\thanksref{label1}}
% \thanks[label1]{}
% \author{Name\corauthref{cor1}\thanksref{label2}}
% \ead{email address}
% \ead[url]{home page}
% \thanks[label2]{}
% \corauth[cor1]{}
% \address{Address\thanksref{label3}}
% \thanks[label3]{}

\title{Some properties of the Lerch family of discrete distributions}\thanks{This work was supported in part by U.S. Public Health Service Grant RO1-GM30054 from the National Institutes of Health.}

% use optional labels to link authors explicitly to addresses:
% \author[label1,label2]{}
% \address[label1]{}
% \address[label2]{}

\author[both,sva]{Sergej V. Aksenov}
\ead{aksenov@freeshell.org}
\and \author[both,mas]{Michael A. Savageau\corauthref{cmas}}
\ead{masavageau@ucdavis.edu}
\corauth[cmas]{Corresponding author}
\address[both]{Department of Microbiology and Immunology, University of Michigan, Ann Arbor, MI 48109}
\address[sva]{Current address: Gene Network Sciences, 31 Dutch Mill Rd., Ithaca, NY 14850}
\address[mas]{Current address: Department of Biomedical Engineering, University of California, One Shields Ave., 1026 Academic Surge, Davis, CA 95616}

\begin{abstract}
% Text of abstract
We extend the definition of the Lerch distribution to the set of nonnegative integers for greater applicability to modeling count data. We express its properties in terms of Lerch's transcendent, and study its unimodality, hazard function and variance-to-mean ratio.

\end{abstract}

\begin{keyword}
% keywords here, in the form: keyword \sep keyword
Discrete distribution \sep Lerch distribution \sep Zipf distribution \sep Zipf-Mandelbrot distribution \sep Lerch's transcendent \sep over-dispersion \sep under-dispersion \sep Mathematica package
% PACS codes here, in the form: \PACS code \sep code

\end{keyword}

\end{frontmatter}

% main text
\section{Introduction}
\label{}
\cite{ZoAl95} introduced a three-parameter, discrete univariate Lerch distribution that is defined on the set of positive integers and is a generalization of related Zipf \citep{Zi49}, Zipf-Mandelbrot \citep{Ma83}, and Good \citep{Go53} distributions. These distributions have been used as models in ecology, linguistics, information science, and statistical physics.  This generalization was made possible by realizing that the probability mass functions (p.m.f.s) $p_x$ of the Zipf, Zipf-Mandelbrot, and Good distributions have the form of the terms in the Riemann zeta ($p_x \sim x^{-s}$), generalized zeta [$p_x \sim (x+v)^{-s}$], and Jonqui\`{e}re's ($p_x \sim z^x x^{-s}$) functions, respectively, which parametrically are special cases of Lerch's transcendent [$p_x \sim z^x (x+v)^{-s}$].  The Lerch distribution is expected to provide a better fit in problems arising in these disciplines because its extra parameters can better capture the shape of the data, which was demonstrated with the zero-truncated Lerch distribution applied to surname data \citep{ZoAl95}.  However, in many other problems the random variable can have value zero, which prompted us to extend the definition of the Lerch distribution to the set of nonnegative integers.  This extension is in fact more natural as a definition, because the series function of the Lerch distribution is Lerch's transcendent, which is defined as an infinite sum with an index running from zero~\citep{MaObSo66}. Truncation of Lerch distribution, including the zero-truncated case studied by \cite{ZoAl95}, can then be easily carried out if needed by manipulating the arguments of Lerch's transcendent (see auxiliary information \citep{Ak02}). Moreover, we show that all functions of the Lerch distribution can be expressed in closed form in terms of Lerch's transcendent. We studied a novel and very effective convergence acceleration technique for computation of Lerch's transcendent elsewhere \citep{AkSa03}; we implemented it in freely available C and Mathematica code \citep{Ak02}. Lerch's transcendent is also available in most major computation software environments including Mathematica \citep{Wo96}. This makes calculations with the Lerch distribution easier in practice.

Motivation for using the Lerch distribution includes the following considerations. \cite{KuTo92} proposed to use the Good distribution to model survival processes because its hazard function can be constant, monotonically decreasing or monotonically increasing depending on the value of one parameter, and dispersal processes because its variance can be greater, equal, or less than the mean.  This allows modeling with a single distribution instead of the combination of three, negative binomial, Poisson, and binomial, as has often been the case in the ecological literature, and instead of compound, generalized Poisson \citep{Co89}, and adjusted generalized Poisson \citep{GuGuTr96}, and other contagious distributions \citep{Ne39,Fe43}.  However, since the Good distribution is not defined for the zero class, it is of no use for certain count problems where the zero class is important.  We show that the general Lerch distribution that includes zero is also flexible with respect to its hazard function and variance-to-mean ratio, and can therefore be successfully applied in those cases. An additional advantage over compound and adjusted distributions is the use of a single distributional form over the whole range of data.  We demonstrate the utility and better goodness-of-fit of the Lerch distribution with an example of under-dispersed data for numbers of sea-urchin sperm in eggs. Calculations for this example are collected in a Mathematica  notebook that is available from the corresponding author.  They have been performed with a specially written Mathematica package LerchDistribution.m that extends the scope of standard statistical functions  of the Mathematica software system to include the Lerch distribution and adds the capability to fit Lerch parameters to data.  The package is freely available for download \citep{Ak02}.

\section{Definition and the distribution functions}
\label{}
The p.m.f. of the Lerch distribution is given by
\begin{equation} \label{lerchpmf}
Pr(X=x) = p_x = \frac{cz^x}{(v+x)^s}
\end{equation}
where $c$ is the normalization constant.  Support of the general Lerch distribution is the set of nonnegative integers $x=0,1,\dots$.  In order for the distribution to be proper, all probabilities $Pr(X=x)$~(\ref{lerchpmf}) have to sum to one:
\begin{equation} \label{proper}
\sum_{x=0}^{\infty} Pr(X=x) = 1
\end{equation}
and so we obtain that the normalization constant is Lerch's transcendent
\begin{equation} \label{lerchnorm}
c^{-1} = \Phi (z,s,v)
\end{equation}
which is defined by the following series \citep{MaObSo66}
\begin{equation} \label{lerchdef}
\Phi(z,s,v) = \sum_{n=0}^{\infty} \frac{z^n}{\left( n+v \right) ^s} \qquad |z| < 1, \quad v \neq 0, -1, \dots
\end{equation}
Also the p.m.f.s~(\ref{lerchpmf}) have to be positive which is ensured first of all by requiring $z>0$.  Likewise, we have to require that $v>0$, because for $v<0$ the p.m.f. alternates in sign for $0\leq x <[|v|]+1$ for integer $s$ and is complex-valued for real $s$, where $[\cdot]$ signifies taking the integer part.  

The cumulative distribution function (c.d.f.) of the Lerch distribution is defined as 
\begin{equation} \label{df}
F(x) = Pr(X \leq x) = 1 - z^{x+1} \frac{\Phi(z,s,v+x+1)}{\Phi(z,s,v)}
\end{equation}
which is obtained using the functional relation  \citep{MaObSo66}:
\begin{equation} \label{lerchrel}
\Phi(z,s,v) = z^m \Phi(z,s,m+v) + \sum_{x=0}^{m-1} \frac{z^x}{(x+v)^s}
\end{equation}
While the c.d.f. of a discrete distribution by definition accepts values only at a countable set of nonnegative integers $x$, we can consider a continuous c.d.f. of the same form as~(\ref{df}) defined on a real line $[0,\infty)$; this continuous c.d.f. of course agrees with our discrete c.d.f. at the nonnegative integers. The definition in~(\ref{df}) implies that this c.d.f., at each integer $x$, is continuous to the right and discontinuous to the left of $x$.  Thus the $q$th quantile of the Lerch distribution is the value $x_q$ such that 
\begin{equation} \label{quantile}
F(x-1)<x_q\leq F(x)
\end{equation}
The solution of this inequality is uniquely defined with probability one, and can be obtained numerically. This also gives an algorithm for random number generation from the Lerch distribution (see the auxiliary information \citep{Ak02}).

The hazard function is defined as 
\begin{equation} \label{haz}
h(x) = \frac{Pr(X=x)}{Pr(X \geq x)} = \frac1{(v+x)^s \Phi(z,s,v+x)}
\end{equation}

The probability generating function (p.g.f.) for $|y| \leq 1$ is
\begin{equation} \label{lerchpgf}
G(y) = \sum_{x=0}^{\infty} y^x Pr(X=x) = \frac{\Phi(yz,s,v)}{\Phi(z,s,v)}
\end{equation}

The moment generating function (m.g.f.) is defined by
\begin{equation} \label{lerchmgf}
M(t) = G(e^t) = \frac{\Phi(ze^t,s,v)}{\Phi(z,s,v)}
\end{equation}
The moments are then obtained as coefficients of the Taylor series expansion of~(\ref{lerchmgf}) in terms of $t$ (see the auxiliary information \citep{Ak02}).

 \section{Some properties}
\label{}

\subsection{Unimodality}
The Lerch distribution, defined in~\textup{(\ref{lerchpmf})}, is strongly unimodal if $s<0$ and $v\geq 1$. These conditions can be derived as follows.

The necessary and sufficient condition for the p.m.f. $p_x$ of any discrete distribution to be strongly unimodal was shown by \cite{KeGe71}:
\begin{equation}
p_x^2\geq p_{x-1} p_{x+1} \label{nec}
\end{equation}
Substituting the p.m.f. of the Lerch distribution~(\ref{lerchpmf}) into the Equation~(\ref{nec}) we obtain
\begin{equation}
\left( 1 - \frac{1}{(v+x)^2} \right)^s \geq 1 \label{nec1}
\end{equation}
Inequality~(\ref{nec1}) can only be satisfied if $s<0$.  Furthermore, the condition
\begin{displaymath}
1 - \frac{1}{(v+x)^2} \geq 0
\end{displaymath}
leads to $v\geq 1$ for all $x=0,1,\dots$, which completes the conditions for strong unimodality.

The mode of the Lerch distribution, defined in~\textup{(\ref{lerchpmf})}, is at $x=[ 1/(z^{1/s} - 1) - v ]+1 = [ 1/(1 - z^{-1/s}) - v ]$, where $[\cdot]$ signifies taking the integer part, provided the conditions of strong unimodality are satisfied. The mode can be obtained by resolving the general inequalities for a discrete distribution with the p.m.f. $p_x$ that has a mode at $x=a$: 
\begin{eqnarray}
p_{x+1}\leq p_x,\qquad \mathrm{all}\, x\geq a \nonumber \\
p_x\geq p_{x-1},\qquad \mathrm{all}\, x\leq a \label{unimod}
\end{eqnarray}
(See auxiliary information \citep{Ak02} for complete derivation.)

\subsection{Monotonicity of the hazard function}
Depending on the sign of parameter $s$, the hazard function can be constant, monotonically decreasing or monotonically increasing.  This can be shown as follows.  First, take the derivative of $h(x)$ with respect to $x$:
\begin{equation} \label{derhaz}
\frac{d h(x)}{d x} = \frac{s}{h(x)} \frac{\Phi(z,s+1,v+x) (v+x) - \Phi(z,s,v+x)}{\Phi(z,s,v+x) (v+x)}
\end{equation}
Now, represent Lerch's transcendents in the numerator of~(\ref{derhaz}) as infinite sums and join these into a single sum, and write the numerator as
\begin{equation} \label{numderhaz}
-z\sum_{k=0}^{\infty} \frac{k z^{k-1}}{(v+x+k)^{s+1}}
\end{equation}
which is evidently negative.  Finally, as $h(x)$ and the denominator in~(\ref{derhaz}) are positive, the sign of the derivative is determined by the sign of $s$: hazard function will be constant if $s=0$, monotonically decreasing if $s>0$, and monotonically increasing if $s<0$.

\subsection{Variance-to-mean ratio}
In order to calculate the variance-to-mean ratio, we need first and second moments. Mean and variance are calculated by differentiating the m.g.f. in~(\ref{lerchmgf}) with respect to $t$ and using relationships between central and uncorrected moments as per \citep{StOr94} (see the auxiliary information \citep{Ak02} for complete derivation). The mean is given by
\begin{equation} \label{mean}
\mu = \mu_1^{\prime} = \frac{\Phi(z,s-1,v)}{\Phi(z,s,v)} - v
\end{equation}
The variance is given by
\begin{equation} \label{var}
\sigma^2 = \mu_2 = (v+\mu)^2 + \frac{\left( -2 (v+\mu) \Phi(z,s-1,v) + \Phi(z,s-2,v) \right)}{\Phi(z,s,v)}
\end{equation}

Now we prove that for a suitable selection of parameters of the Lerch distribution its variance can be greater, equal or less than its mean, thus making it useful for analysis of over- and under-dispersed data.  Consider the Taylor series expansion of the mean and the variance of the Lerch distribution given by Equations~(\ref{mean}) and~(\ref{var}), respectively, around $z=0$ to order $z^1$ and $z^2$, respectively:
\begin{eqnarray} \label{mvtaylor}
\mu & = & z \frac{v^s}{(1+v)^s} + O(z^2) \nonumber \\
\sigma^2 & = & z \frac{v^s}{(1+v)^s} + z^2 2 (v+1) \left( 2\frac{v^s}{(2+v)^s} - \frac{v^{2s}}{(1+v)^{2s}} \right) + O(z^3)
\end{eqnarray}
Now calculate the ratio of variance to mean to order $z^1$ using Equations~(\ref{mvtaylor}):
\begin{equation} \label{vmratio}
\frac{\sigma^2}{\mu} = 1+z 2 (v+1) \left[ 2 \left(\frac{1+v}{2+v}\right)^s-\left(\frac{v}{1+v}\right)^s \right]
\end{equation}
Clearly, the ratio of variance to mean will be less, equal or greater than unity, depending on the sign of the expression in square brackets in Equation~(\ref{vmratio}).  This expression can change its sign only if $s<0$.  Thus, we obtain that for any $v>0$ and $z>0$, which is small in Taylor's sense, there is $s<0$ such that
\begin{equation} \label{vmratcond}
\frac{\sigma^2}{\mu} <,=,> 1 \qquad iff \qquad s<,=,> -\frac{\log 2}{\log \left(1+\frac1{v^2+2v}\right)}
\end{equation}
Furthermore, based on numerical experiments, we conjecture that for $s > -\log 2/\log (1+1/(v^2+2v))$, the variance-to-mean ratio will remain greater than unity for larger deviations of $z$ from zero.  However, for $s < -\log 2/\log (1+1/(v^2+2v))$, the ratio remains less than one only for small deviations of $z$ from 0; for larger deviations in $z$ the ratio becomes greater than unity.

\section{Example}
\label{}

Consider fertilization of sea-urchin eggs by sperm (rows 1 and 4  in Table~\ref{tab:urchin}) \citep{Mo75}, which is an example of underdispersed data that is fitted well by the Lerch distribution but not fitted at all by the generalized Poisson distribution.  Experiments are done by examining the number of fertilized eggs in a batch at two different times, 40 and 180 sec. Poisson model failed to fit the data, suggesting a nonrandom pattern of fertilization. \cite{JaKeSc79} fit the generalized Poisson disrtribution that has the p.f. \citep{Co89}
\begin{equation} \label{genp}
p_x = \left\{ 
\begin{array}{ll}
\frac{\theta (\theta + x \lambda)^{x-1} \mathrm{e}^{-\theta -x \lambda}}{x!}, & x=0,1,2,\dots \\
0,& x>\left[-\frac{\theta}{\lambda}\right] \quad \mathrm{when} \quad \lambda<0
\end{array} \right.
\end{equation}
where $[\cdot]$ signifies the integer part operation, with parameters $\theta=1.0077$ and $\lambda=-0.3216$ for 40 sec. and $\theta=2.4654$ and $\lambda=-1.0893$ for 180 sec. (see rows 2 and 5 in Table~\ref{tab:urchin}, respectively, for calculated frequencies).  However, since $\lambda<0$ for both data sets, the generalized Poisson distribution should be truncated at 3 and 2, respectively.  Clearly the truncation points being less than the maximum observed class indicate that generalized Poisson distribution fails to fit the whole range of data.  We fit the Lerch distribution by minimizing the Pearson $X^2$ statistic.  The best-fit parameters are $z=0.00773867$, $s=-8.26894$, $v=1.11633$ for 40 sec. and $z=0.0835808$, $s=-1.15174$, $v=0.00468234$ for 180 sec. (see rows 3 and 6 in Table~\ref{tab:urchin} for calculated frequencies).  Since there are only four nonzero classes in both data sets, the Chi-squared distribution would have zero d.f., and thus we use the sum of squared deviations as a goodness-of-fit measure.  For 40 sec. the sum is 0.000372774 and for 180 sec. it is 0.000162184. 

The above example demonstrates that the Lerch distribution can provide a useful model with better goodness-of-fit for over- and under-dispersed data (see the auxiliary information  \citep{Ak02} for additional examples), for which modifications of the Poisson distribution or the reduced Lerch (Good) model perform less well or even fail, as well as for some rank-abundance ecological data, for which the reduced Lerch (Zipf-Mandelbrot) model does not fit.  The ultimate reason is that the Lerch distribution is a generalization of the Good and Zipf-Mandelbrot distributions (which were originally proposed as models for these problems) and is endowed with greater flexibility to accomodate the odd distributional shapes often observed in real problems.  A comment on the method of parameter estimation is in order.  We found in practice that the most convenient and accurate method of parameter estimation is minimization of the Pearson $X^2$ statistic, which works well for data even with a few frequency classes.  The moment and maximum likelihood methods derived in the auxiliary information  \citep{Ak02} are more difficult to use in such situations because of higher variance of sample moments and expectations in a small sample setting.  An added difficulty is the slow convergence of moment and maximum likelihood equations that involve Lerch's transcendent.  These difficulties are not specific to the Lerch distribution and are in fact common for other distributions that are based on special functions.  However, we found that these methods perform well with more rich data sets (not shown).

% The Appendices part is started with the command \appendix;
% appendix sections are then done as normal sections
% \appendix

% \section{}
% \label{}

% Bibliographic references with the natbib package:
% Parenthetical: \citep{Bai92} produces (Bailyn 1992).
% Textual: \citet{Bai95} produces Bailyn et al. (1995).
% An affix and part of a reference:
%   \citep[e.g.][Ch. 2]{Bar76}
%   produces (e.g. Barnes et al. 1976, Ch. 2).

%\begin{thebibliography}{}

% \bibitem[Names(Year)]{label} or \bibitem[Names(Year)Long names]{label}.
% (\harvarditem{Name}{Year}{label} is also supported.)
% Text of bibliographic item

%\bibitem[]{}

%\end{thebibliography}

\bibliography{REFS}

\begin{table} 
\begin{center}
\begin{tabular}{lllllll}
\hline
Col. & Counts & 0 & 1 & 2 & 3 & 4  \\
\hline
\multicolumn{6}{l}{40 sec.}\\
\hline
1& Observed & 28 & 44 & 7 & 1 & 0 \\
2& Gen. Poisson & 29.2046 & 40.5931 & 10.2044 & 0.0236894 & -- \\
3& Lerch & 28.0876 & 43.0737 & 8.17627 & 0.631919 & 0.0295335 \\
\hline
\multicolumn{6}{l}{180 sec.}\\
\hline
4& Observed & 2 & 81 & 15 & 1 & 1 \\
5& Gen. Poisson & 8.49748 & 62.2665 & 26.5388 & -- & -- \\
6& Lerch & 1.99482 & 80.7898 & 14.9625 & 1.99311 & 0.23192 \\
\hline
\end{tabular}
\end{center}
\vspace{0.5in}
\caption{Distribution of sea-urchin sperm on eggs at two different times of exposure.  The observations are frequencies of sperm on eggs, frequency classes are numbers of sperm found on an egg.  Observed frequencies come from \citep{Mo75}. Predicted frequencies come from the best-fit generalized Poisson and Lerch distributions.  See text for details.} \label{tab:urchin}
\end{table}

\end{document}

\begin{frontmatter}

% Title, authors and addresses

% use the thanksref command within \title, \author or \address for footnotes;
% use the corauthref command within \author for corresponding author footnotes;
% use the ead command for the email address,
% and the form \ead[url] for the home page:
% \title{Title\thanksref{label1}}
% \thanks[label1]{}
% \author{Name\corauthref{cor1}\thanksref{label2}}
% \ead{email address}
% \ead[url]{home page}
% \thanks[label2]{}
% \corauth[cor1]{}
% \address{Address\thanksref{label3}}
% \thanks[label3]{}

\title{Auxiliary information for article Some properties of the Lerch family of discrete distributions}\thanks{This work was supported in part by U.S. Public Health Service Grant RO1-GM30054 from the National Institutes of Health.}

% use optional labels to link authors explicitly to addresses:
% \author[label1,label2]{}
% \address[label1]{}
% \address[label2]{}

\author[both,sva]{Sergej V. Aksenov}
\ead{aksenov@freeshell.org}
\and \author[both,mas]{Michael A. Savageau\corauthref{cmas}}
\ead{masavageau@ucdavis.edu}
\corauth[cmas]{Corresponding author}
\address[both]{Department of Microbiology and Immunology, University of Michigan, Ann Arbor, MI 48109}
\address[sva]{Current address: Gene Network Sciences, 31 Dutch Mill Rd., Ithaca, NY 14850}
\address[mas]{Current address: Department of Biomedical Engineering, University of California, One Shields Ave., 1026 Academic Surge, Davis, CA 95616}

\begin{abstract}
% Text of abstract
This text contains derivation of expression for the mode of the Lerch distribution, algorithm for random number generation, calculation of moments of the Lerch distribution, properties of truncated Lerch distributions, formulas for parameter estimation by method of moments and maximum likelihood. Finally, we give some examples showing superior goodness-of-fit when fitting Lerch distribution to various data.
\end{abstract}

\begin{keyword}
% keywords here, in the form: keyword \sep keyword
Discrete distribution \sep Lerch distribution \sep Zipf distribution \sep Zipf-Mandelbrot distribution \sep Lerch's transcendent \sep   over-dispersion \sep under-dispersion \sep Mathematica package
% PACS codes here, in the form: \PACS code \sep code

\end{keyword}

\end{frontmatter}

% main text
\section{Mode}

The p.m.f. of the Lerch distribution is given by
\begin{equation} \label{lerchpmf}
Pr(X=x) = p_x = \frac{cz^x}{(v+x)^s}
\end{equation}
where $c$ is Lerch's transcendent
\begin{equation} \label{lerchnorm}
c^{-1} = \Phi (z,s,v)
\end{equation}
which is defined by the following series \citep{MaObSo66}
\begin{equation} \label{lerchdef}
\Phi(z,s,v) = \sum_{n=0}^{\infty} \frac{z^n}{\left( n+v \right) ^s} \qquad |z| < 1, \quad v \neq 0, -1, \dots
\end{equation}

\begin{prop} 

The mode of the Lerch distribution, defined in~\textup{(\ref{lerchpmf})}, is at $x=[ 1/(z^{1/s} - 1) - v ]+1 = [ 1/(1 - z^{-1/s}) - v ]$, where $[\cdot]$ signifies taking the integer part, provided the conditions for strong unimodality as per main text are satisfied.

\end{prop}
%
%\begin{proof}

\noindent {\bf Proof.} A discrete distribution with the p.m.f. $p_x$ has a mode at $x=a$ if the following inequalities hold: 
\begin{eqnarray}
p_{x+1}\leq p_x,\qquad \mathrm{all}\, x\geq a \nonumber \\
p_x\geq p_{x-1},\qquad \mathrm{all}\, x\leq a \label{unimod}
\end{eqnarray}
Using Equation~(\ref{lerchpmf}) we rewrite these inequalities as
\begin{equation}
z \left( \frac{v+x}{v+x+1} \right)^s \leq 1, \qquad z \left( \frac{v+x-1}{v+x} \right)^s \geq 1 \label{dermod}
\end{equation}
Provided that the conditions for unimodality are satisfied, these inequalities can be resolved to obtain the mode.

First, these inequalities can be rearranged as (we use the fact that raising quantities that are less than one to negative power reverses the inequality sign)
\begin{displaymath}
\frac1{v+x} \leq z^{1/s} - 1,\qquad \frac1{v+x} \geq 1 - z^{-1/s}
\end{displaymath}
Considering that $0<z<1$, $s<0$ and $v\geq 1$, we see that both the left and right hand sides of both inequalities are positive.  Rearranging we obtain further
\begin{displaymath}
x \geq \frac1{z^{1/s} - 1} - v,\qquad x \leq \frac1{1 - z^{-1/s}} - v
\end{displaymath}
We can see that the (real) endpoints of the half-intervals defined by these inequalities bracket the mode which is a unique integer, because their difference is exactly one:
\begin{displaymath}
\frac1{1 - z^{-1/s}} - \frac1{z^{1/s} - 1} = \frac{z^{1/s}}{z^{1/s} - 1} - \frac1{z^{1/s} - 1} = 1
\end{displaymath}
Taking the integer part of either of the endpoints we obtain the expressions for the mode.$\Box$
%\end{proof}

\section{Random number generation}

Random number generation from a discrete distribution can in principle be accomplished with a number of techniques \citep{De86}.  We propose a direct way to obtain the discrete Lerch random variate using inversion by correction.  In the inversion by correction method, one starts with a c.d.f. $H(x)$ that is "close" to the true c.d.f $F(x)$
\begin{equation} \label{df}
F(x) = 1 - z^{x+1} \frac{\Phi(z,s,v+x+1)}{\Phi(z,s,v)}
\end{equation}
and there is an "easy" algorithm to generate a random number from $H(x)$ (e.g., straightforward inversion).  Given $H(x)$, one  samples from $H(x)$ and then corrects by sequential search to obtain a random variate from $F(x)$ \citep{De86}.  Obtaining  a first estimate of $x$ in this way serves to shorten the ordinary sequential search.  Additional savings of computational time can be obtained if $H(x)$ can be constructed in such a way that $F(x)$ and $H(x)$ are stochastically ordered, i.e. $H(x)\leq F(x)$ or $H(x)\geq F(x)$ for all $x$.  In this case we know beforehand if we have to search up or down.  In the case of the Lerch c.d.f., it is possible to construct a stochastically ordered $H(x)$ as: 
\begin{equation} \label{stochh}
H(x)=1-z^{x+1} \frac{\Phi(z,s,v+1)}{\Phi(z,s,v)}
\end{equation}
$H(x)$ and $F(x)$ are stochastically ordered because $\Phi(z,s,v+1) > \Phi(z,s,v+x+1)$ and thus $H(x)\leq F(x)$ if $s>0$, and because $\Phi(z,s,v+1) \leq \Phi(z,s,v+x+1)$ and thus $H(x)\geq F(x)$ if $s\leq0$, for all $x$ as $x\rightarrow \infty$.

The algorithm to generate a random Lerch variate is then as follows.  
\begin{algorithm} \label{rng}

Inversion by correction: Lerch distribution($z,s,v$)

\begin{enumerate}
\item Generate a uniform $(0,1)$ random variate $U$.  
\item Calculate a first guess $x$ by inversion and truncation of a continuous $H(x)$,
\begin{displaymath}
x=\left[1+\frac{\log (1-U) \Phi(z,s,v)/\Phi(z,s,v+1)}{\log z} \right]-1
\end{displaymath}
where $[\cdot]$ signifies the integer part operation.  
\item Finally, if $s>0$, $x$ is adjusted down sequentially, $x=x-1$, until $U> F(x-1)$.  If $s\leq 0$, we adjust up sequentially, $x=x+1$, until $U\leq F(x)$.  
\end{enumerate}

\end{algorithm}

The time it takes to generate a random number is on average proportional to the expected number of calculations of $F$, $E(C)$, which in turn is related to the number of comparisons during the search and hence to the expectation $E(x)$, $E(C)=1+E(x)$ \citep{De86}.

\begin{thm}

The expected number of computations, $E(C)$, for the sequential search in Algorithm~\ref{rng} is less than the expected number in the ordinary sequential search, for all parameter values where the Lerch distribution is defined.

\end{thm}

%\begin{proof}

\noindent {\bf Proof.} For the ordinary sequential search the expected number of computations is 
\begin{displaymath}
E_1 = 1+E(x) = 1+\sum_x (1-F(x)) = 1+\frac{\sum_{x=0}^{\infty} z^{x+1} \Phi(z,s,v+x+1)}{\Phi(z,s,v)}
\end{displaymath}
and for the inversion by correction method it is
\begin{eqnarray}
E_2 & = & 1+E(x) = 1+\sum_x |F(x)-H(x)| = 1+ \nonumber \\
&& \frac{\sum_{x=0}^{\infty} z^{x+1} |\Phi(z,s,v+1)-\Phi(z,s,v+x+1)|}{\Phi(z,s,v)} \nonumber
\end{eqnarray}

If $H\geq F$, then
\begin{equation}
E_1-E_2 = \frac{\Phi(z,s,v+1)}{\Phi(z,s,v)} \sum_{x=0}^{\infty} z^{x+1} = \frac{\Phi(z,s,v+1)}{\Phi(z,s,v)} \frac{z}{1-z} \label{hgf}
\end{equation}
which is always a positive quantity, and thus $E_1>E_2$.

If $H \leq F$, then
\begin{equation}
E_1-E_2 = \sum_{x=0}^{\infty} (1 - (2F(x)-H(x))) = \sum_{x=0}^{\infty} (1-G(x)) \label{hlf}
\end{equation}
which is an expectation of $x$ with respect to the c.d.f. $G(x)$ on a set of nonnegative integers, and thus a positive quantity.  Hence $E_1>E_2$.

Here $G(x)=2F(x)-H(x)$ is a valid c.d.f. because $G(x)$ is a nondecreasing function of $x$, $G(0)=0$ and $\lim_{x\rightarrow \infty} G(x) = 1$.$\Box$
%\end{proof}

\begin{cor} \label{cor1}

The difference between the expected numbers of computations, $E(C)$, for the ordinary sequential search and the inversion by correction Algorithm~\ref{rng} is less if $s>0$ than if $s\leq 0$.

\end{cor}

%\begin{proof}

\noindent {\bf Proof.} By using Equations~(\ref{df}) and~(\ref{stochh}) we rewrite Equation~(\ref{hlf}), which corresponds to $s>0$, as
\begin{eqnarray}
E_1-E_2 & = & 2 \frac{\sum_{x=0}^{\infty} z^{x+1} \Phi(z,s,v+x+1)}{\Phi(z,s,v)} - \frac{\Phi(z,s,v+1)}{\Phi(z,s,v)} \sum_{x=0}^{\infty} z^{x+1} < \nonumber \\
&& \frac{\Phi(z,s,v+1)}{\Phi(z,s,v)} \frac{z}{1-z} \label{corl}
\end{eqnarray}
The inequality is obtained by using the fact that for $s>0$ we have $\Phi(z,s,v+1) > \Phi(z,s,v+x+1)$.

Finally, observe that the difference in Equation~(\ref{corl}) is less than the difference in Equation~(\ref{hgf}).$\Box$
%\end{proof}

Corollary~\ref{cor1} implies that the gain in efficiency by using inversion by correction vs. ordinary sequential search is diminished for $s>0$; the gain is higher if $s\leq 0$.  Obviously, Algorithm~\ref{rng} can also be used for calculation of a $q$th quantile replacing $U$ with $q$.

\section{Moments}

The moment generating function (m.g.f.) is defined by
\begin{equation} \label{lerchmgf}
M(t) = G(e^t) = \frac{\Phi(ze^t,s,v)}{\Phi(z,s,v)}
\end{equation}

The factorial (descending) moment generating function (f.m.g.f.) is defined by
\begin{equation} \label{lerchfmgf}
G(1+t) = \frac{\Phi(z(t+1),s,v)}{\Phi(z,s,v)}
\end{equation}
and the factorial (descending) moments are the coefficients of the Taylor series of~(\ref{lerchfmgf}).

The $r$th uncorrected moment can be obtained by differentiating $r$ times the m.g.f. in~(\ref{lerchmgf}) with respect to $t$ and letting $t=0$:
\begin{eqnarray} \label{defuncor}
\mu_r^{\prime} & = & \left( \frac{d^r G(e^t)}{dt^r} \right)_{t=0} = \frac1{\Phi(z,s,v)} \left( \frac{d^r \Phi(e^t z,s,v)}{dt^r} \right)_{t=0} \nonumber \\ 
& & = \frac1{\Phi(z,s,v)} \sum_{x=0}^{\infty} \frac{x^r z^x}{(x+v)^s}, \qquad r=1,2,\dots
\end{eqnarray}
By using the binomial expansion for $x^r = (x+v-v)^r = \sum_{j=0}^{r} {r \choose j} (x+v)^{j} (-v)^{r-j}$ and changing the order of summation \citep{StOr94} we obtain further 
\begin{equation} \label{compuncor}
\mu_r^{\prime} = \frac1{\Phi(z,s,v)} \sum_{j=0}^{r} {r \choose j} (-v)^{r-j} \Phi(z,s-j,v)
\end{equation}
Thus, the mean is given by
\begin{equation} \label{mean}
\mu = \mu_1^{\prime} = \frac{\Phi(z,s-1,v)}{\Phi(z,s,v)} - v
\end{equation}
Similarly, the second and third uncorrected moments are given by
\begin{eqnarray} \label{uncor23}
\mu_2^{\prime} & = & v^2 -2 v \frac{\Phi(z,s-1,v)}{\Phi(z,s,v)} + \frac{\Phi(z,s-2,v)}{\Phi(z,s,v)}  \nonumber \\
\mu_3^{\prime} & = & -v^3 + 3 v^2 \frac{\Phi(z,s-1,v)}{\Phi(z,s,v)} \nonumber \\
&& - 3 v \frac{\Phi(z,s-2,v)}{\Phi(z,s,v)} + \frac{\Phi(z,s-3,v)}{\Phi(z,s,v)}
\end{eqnarray}

The central moments are obtained from the uncorrected moments using the following relationship \citep{StOr94}:
\begin{equation} \label{uncortocentral}
\mu_r = \sum_{k=0}^{r} (-1)^k {r \choose k} \mu_{r-k}^{\prime} (\mu_1^{\prime})^k
\end{equation}
Substituting the uncorrected moments $\mu_{r-k}^{\prime}$ from Equation~(\ref{compuncor}) into Equation~(\ref{uncortocentral}) and rearranging terms we obtain
\begin{equation} \label{central1}
\mu_r = \sum_{k=0}^{r} \sum_{j=0}^{r-k} {r \choose k} {r-k \choose j} (-v)^{r-k-j} (-\mu_1^{\prime})^k \frac{\Phi(z,s-j,v)}{\Phi(z,s,v)}
\end{equation}
Now exchange the order of the summations, use the relationship 
\begin{displaymath}
{r\choose k} {r-k\choose j} = {r-j\choose k} {r\choose j}
\end{displaymath}
and then collapse the resulting binomial expansion to obtain  
\begin{equation} \label{central20}
\mu_r = \frac1{\Phi(z,s,v)} \sum_{j=0}^{r} (-1)^{r-j} {r\choose j} \Phi(z,s-j,v) (v+\mu_1^{\prime})^{r-j}
\end{equation}
Finally, substituting the mean~(\ref{mean}), we obtain
\begin{equation} \label{central2}
\mu_r = \frac1{\Phi(z,s,v)} \sum_{j=0}^{r} (-1)^{r-j} {r\choose j} \Phi(z,s-j,v) \left( \frac{\Phi(z,s-1,v)}{\Phi(z,s,v)} \right)^{r-j}
\end{equation}
Thus, the variance is given by
\begin{equation} \label{var}
\sigma^2 = \mu_2 = (v+\mu)^2 + \frac{\left( -2 (v+\mu) \Phi(z,s-1,v) + \Phi(z,s-2,v) \right)}{\Phi(z,s,v)}
\end{equation}
Similarly, the third and fourth central moments are
\begin{eqnarray} \label{central34}
\mu_3 & = & -(v+\mu)^3 + \frac1{{\Phi(z,s,v)}} \left( 3 (v+\mu)^2 \Phi(z,s-1,v) - \right. \nonumber \\
&& \left. 3 (v+\mu) \Phi(z,s-2,v) + \Phi(z,s-3,v) \right) \nonumber \\
\mu_4 & = & (v+\mu)^4 + \frac1{{\Phi(z,s,v)}} \left( -4 (v+\mu)^3 \Phi(z,s-1,v) + \right. \nonumber \\
&& \left. 6 (v+\mu)^2 \Phi(z,s-2,v) - 4 (v+\mu) \Phi(z,s-3,v) + \right. \nonumber \\
&& \left. \Phi(z,s-4,v) \right) 
\end{eqnarray}

The ratios of some central moments are commonly used as indices for the shape of the distribution.  Namely, skewness is defined as
\begin{equation}\label{skew}
\alpha_3 = \frac{\mu_3}{\mu_2^{3/2}}
\end{equation}
Kurtosis is defined as
\begin{equation} \label{kurt}
\alpha_4 = \frac{\mu_4}{\mu_2^2}
\end{equation}
Skewness and kurtosis can be calculated by substituting the mean and the central moments from Equations~(\ref{mean}), (\ref{var}) and~(\ref{central34}), respectively.

The $r$th factorial (descending) moments can be obtained from the f.m.g.f. by differentiation $r$ times with respect to $t$ and letting $t=0$:
\begin{eqnarray} \label{factmom}
\mu_{[r]}^{\prime} & = & \left( \frac{d^r G(1+t)}{dt^r} \right)_{t=0} = \frac1{\Phi(z,s,v)} \left( \frac{d^r \Phi((1+t) z,s,v)}{dt^r} \right)_{t=0} \nonumber \\ 
& & = \frac1{\Phi(z,s,v)} \sum_{x=0}^{\infty} \frac{x^{(r)} z^x}{(x+v)^s}, \qquad r=1,2,\dots \end{eqnarray}
where $x^{(r)}=x(x-1)\dots(x-r+1)$ is the descending factorial.  By inverting Newton's difference formula for the power function $x^r$, the descending factorials can be expressed as functions of powers of $x$ \citep{JoKoKe92}:
\begin{equation} \label{factstir}
x^{(r)} = \sum_{i=0}^r s(r,i) x^i
\end{equation}
where $s(r,i)$ are the Stirling numbers of the first kind \citep{GoLeNeZu76}.  Substituting Equation~(\ref{factstir}) into Equation~(\ref{factmom}) and changing the order of summation we obtain
\begin{equation} \label{fact}
\mu_{[r]}^{\prime} = \sum_{i=0}^r s(r,i) \sum_{x=0}^{\infty} \frac{x^i z^x}{(x+v)^s} = \sum_{i=0}^r s(r,i) \mu_i^{\prime}, \qquad r=1,2,\dots 
\end{equation}
which relates factorial moments to the uncorrected ones.  Thus, the first three factorial (descending) moments are
\begin{eqnarray} \label{fact123}
\mu_{[1]}^{\prime} & = & \mu_1^{\prime} = \frac{\Phi(z,s-1,v)}{\Phi(z,s,v)} - v  \nonumber \\
\mu_{[2]}^{\prime} & = & \mu_2^{\prime} - \mu_1^{\prime} = v(v+1) - \frac{\Phi(z,s-1,v)}{\Phi(z,s,v)} (2v+1) + \frac{\Phi(z,s-2,v)}{\Phi(z,s,v)} \nonumber \\
\mu_{[3]}^{\prime} & = & \mu_3^{\prime}-3 \mu_2^{\prime}+2\mu_1^{\prime} = -v(v+1)(v+2)+ \nonumber \\
&& \frac{\Phi(z,s-1,v)}{\Phi(z,s,v)} (3v^2+6v+2) -3 \frac{\Phi(z,s-2,v)}{\Phi(z,s,v)} (v+1)+ \nonumber \\
&& \frac{\Phi(z,s-3,v)}{\Phi(z,s,v)}
\end{eqnarray}

\section{Truncated Lerch distribution}

The Lerch distribution can be considered in singly or doubly truncated forms, which will change the summation limits in the normalization equation
\begin{equation} \label{proper}
\sum_{x=0}^{\infty} Pr(X=x) = 1
\end{equation}
To calculate the normalization constant in these cases, we make use of the functional relation \citep{MaObSo66}:
\begin{equation} \label{lerchrel}
\Phi(z,s,v) = z^m \Phi(z,s,m+v) + \sum_{x=0}^{m-1} \frac{z^x}{(x+v)^s}
\end{equation}
If the truncation points are $a \geq 0$ and $a \leq b \leq \infty$, the normalization constant is
\begin{equation} \label{lerchnormtrunc}
c^{-1} = \sum_{x=a}^b \frac{z^x}{(v+x)^s} = z^a \Phi(z,s,v+a) - z^{b+1} \Phi(z,s,v+b+1)
\end{equation}
The requirement $v>0$ for the distribution to be proper is replaced with $v>-a$ for a truncated form of the distribution.

For example, in the zero-truncated case we have $a=1$ and $b=\infty$ and $c^{-1} = z \Phi(z,s,v+1)$.  This expression has to be contrasted with the one immediately following Equation~(2.1b) in \citep{ZoAl95} and all other equations in which they use the symbol $\Phi$.  One is cautioned that in \citep{ZoAl95}, the symbol $\Phi$ is referred to as Lerch's transcendent (see Equation~(1.1) in \citep{ZoAl95}), but it lacks the first term of the infinite series and is thus \emph{not} the correct definition of Lerch's transcendent according to Equation~(\ref{lerchdef}).  However, the calculations in \citep{ZoAl95} will be correct if one uses direct summation instead of Lerch's transcendent $\Phi$ as implemented in various software systems.

The c.d.f. is
\begin{equation} \label{dftrunc}
F(x) = c \sum_{n=a}^x \frac{z^n}{(n+v)^s} = \frac{z^a \Phi(z,s,v+a) - z^{x+1} \Phi(z,s,v+x+1)}{z^a \Phi(z,s,v+a) - z^{b+1} \Phi(z,s,v+b+1)}
\end{equation}
where we used Equation~(\ref{lerchnormtrunc}) for the normalization constant.  

The survival function is
\begin{equation} \label{survtrunc}
S(x) = z \frac{z^x \Phi(z,s,v+x+1) - z^b \Phi(z,s,v+b+1)}{z^a \Phi(z,s,v+a) - z^{b+1} \Phi(z,s,v+b+1)}
\end{equation}

The hazard function is
\begin{equation} \label{haztrunc}
h(x) = \left( z (v+x)^s ( \Phi(z,s,v+x+1)-z^{b-x} \Phi(z,s,v+b+1))\right)^{-1}
\end{equation}

The p.g.f. is given by
\begin{equation} \label{pgftrunc}
G(y) = \frac{(yz)^a \Phi(yz,s,v+a) - (yz)^{b+1} \Phi(yz,s,v+b+1)}{z^a \Phi(z,s,v+a) - z^{b+1} \Phi(z,s,v+b+1)}
\end{equation}
For the zero-truncated case we obtain
\begin{equation} \label{pgftrunc1}
G(y) = \frac{y \Phi(yz,s,v+1)}{\Phi(z,s,v+1)}
\end{equation}
which should be contrasted with Equation~(2.2) in \citep{ZoAl95}, where again their $\Phi$ is \emph{not} the correct Lerch's transcendent of Equation~(\ref{lerchdef}).

Uncorrected moments are given by
\begin{eqnarray} \label{compuncortrunc}
\mu_r^{\prime} & = & \frac1{z^a \Phi(z,s,v+a) - z^{b+1} \Phi(z,s,v+b+1)} \times \nonumber \\
&& \sum_{j=0}^{r} {r \choose j} (-v)^{r-j} (z^a \Phi(z,s-j,v+a) - \nonumber \\
&& z^{b+1} \Phi(z,s-j,v+b+1))
\end{eqnarray}
and central moments are given by
\begin{eqnarray} \label{centraltrunc}
\mu_r & = & \frac1{z^a \Phi(z,s,v+a) - z^{b+1} \Phi(z,s,v+b+1)} \times \nonumber \\
&& \sum_{j=0}^{r} (-1)^{r-j} {r\choose j} (z^a \Phi(z,s-j,v+a) - z^{b+1} \Phi(z,s-j,v+b+1)) \times \nonumber \\
&& \left( \frac{z^a \Phi(z,s-1,v+a) - z^{b+1} \Phi(z,s-1,v+b+1)}{z^a \Phi(z,s,v+a) - z^{b+1} \Phi(z,s,v+b+1)} \right)^{r-j}
\end{eqnarray}
For example, in the zero-truncated case one obtains uncorrected moments as
\begin{equation} \label{unc1trunc}
\mu_r^{\prime} = \frac1{\Phi(z,s,v+1)} \sum_{j=0}^{r} {r \choose j} (-v)^{r-j} \Phi(z,s-j,v+1)
\end{equation}
and central moments as
\begin{eqnarray} \label{cen1trunc}
\mu_r & = & \frac1{\Phi(z,s,v+1)} \times \nonumber \\
&& \sum_{j=0}^{r} (-1)^{r-j} {r\choose j} \Phi(z,s-j,v+1) \left( \frac{\Phi(z,s-1,v+1)}{\Phi(z,s,v+1)} \right)^{r-j}
\end{eqnarray}
Equation~(\ref{cen1trunc}) has to be contrasted with Equation~(3.1) in \citep{ZoAl95}, where the running index $j$ incorrectly starts from 1 and $\Phi$ is again \emph{not} the correct Lerch's transcendent of Equation~(\ref{lerchdef}).

\section{Parameter estimation}

Suppose we are given a set of nonnegtaive integer-valued data
\begin{equation} \label{data1}
\{x_i\}, \qquad i=1,\dots,n
\end{equation}
Sample data may also be given as a vector of empirical frequences, 
\begin{equation} \label{data2}
\{f_j=\#(x=x_j)/n\},\qquad x_j=0,1,\dots
\end{equation}
where $\#$ signifies the number of occurences.  Here we consider two methods for fitting the Lerch distribution to data given as~(\ref{data1}) or~(\ref{data2}), the moment and maximum likelihood methods. We implemented both methods in a Mathematica package LerchDistribution.m which is available for download \citep{Ak02url}.

\subsection{Method of moments}

The moment method (MM) for the Lerch distribution involves equating the first three uncorrected moments~(\ref{mean}) and~(\ref{uncor23}) [or central moments~(\ref{var}) and~(\ref{central34}), or factorial moments~(\ref{fact123})] to their sample counterparts 
\begin{eqnarray} \label{samplemom1}
m_r^{\prime} & = & n^{-1} \sum_{i=1}^{n} x_i^r = \sum_j f_j x_j^r \nonumber \\
m_r & = & n^{-1} \sum_{i=1}^{n} (x_i-m_1^{\prime})^r = \sum_j f_j (x_j-m_1^{\prime})^r \nonumber \\
m_{[r]}^{\prime} & = & n^{-1} \sum_{i=1}^{n} x_i^{[r]} = \sum_j f_j x_j^{[r]}
\end{eqnarray}
for uncorrected, central and factorial moments, respectively. 

For example, equating uncorrected moments gives the following set of equations for estimators $\bar{z}$, $\bar{s}$ and $\bar{v}$ (equating central or factorial moments can be done in a similar way):
\begin{eqnarray} \label{uncorest}
m_1^{\prime} & = & \frac{\Phi(\bar{z},\bar{s}-1,\bar{v})}{\Phi(\bar{z},\bar{s},\bar{v})} - \bar{v} \nonumber \\
m_2^{\prime} & = &\bar{v}^2 -2 \bar{v} \frac{\Phi(\bar{z},\bar{s}-1,\bar{v})}{\Phi(\bar{z},\bar{s},\bar{v})} + \frac{\Phi(\bar{z},\bar{s}-2,\bar{v})}{\Phi(\bar{z},\bar{s},\bar{v})} \nonumber \\
m_3^{\prime} & = & -\bar{v}^3 + 3 \bar{v}^2 \frac{\Phi(\bar{z},\bar{s}-1,\bar{v})}{\Phi(\bar{z},\bar{s},\bar{v})} - 3 \bar{v} \frac{\Phi(\bar{z},\bar{s}-2,\bar{v})}{\Phi(\bar{z},\bar{s},\bar{v})} + \nonumber \\
&& \frac{\Phi(\bar{z},\bar{s}-3,\bar{v})}{\Phi(\bar{z},\bar{s},\bar{v})}
\end{eqnarray}
Upon substituting values of sample moments, Equations~(\ref{uncorest}) can be solved numerically to obtain MM estimates $\bar{z}$, $\bar{s}$ and $\bar{v}$.

The approximate asymptotic variance-covariance matrix for MM estimates can be obtained using the delta method \citep{StOr94}, in which one approximates the right-hand sides of Equations~(\ref{uncorest}), $\mu_r^{\prime}(\theta)$ where $\theta=(z,s,v)$, by the first terms of Taylor series expansion.  Then after applying the variance operator we have:
\begin{equation} \label{varmu}
var (\mu_r^{\prime}) = \sum_{j=1}^3 \left( \frac{\partial \mu_r^{\prime}}{\partial \theta_j} \right)^2 var (\theta_j) + \sum_{k=1}^3 \sum_{l=1,l\neq k}^3 \frac{\partial \mu_r^{\prime}}{\partial \theta_k} \frac{\partial \mu_r^{\prime}}{\partial \theta_l} cov (\theta_k,\theta_l)
\end{equation}
Similarly, the covariance of population moments $\mu_r^{\prime}(\theta)$ is
\begin{equation} \label{covmu}
cov (\mu_r^{\prime}, \mu_q^{\prime}) = \sum_{j=1}^3 \frac{\partial \mu_r^{\prime}}{\partial \theta_j} \frac{\partial \mu_q^{\prime}}{\partial \theta_j} var (\theta_j) + \sum_{k=1}^3 \sum_{l=1,l\neq k}^3 \frac{\partial \mu_r^{\prime}}{\partial \theta_k} \frac{\partial \mu_q^{\prime}}{\partial \theta_l} cov (\theta_k,\theta_l)
\end{equation}
To obtain equations for the variances and covariances of $\theta$ we now calculate the variance and covariance of sample moments $m_r^{\prime}$ \citep{StOr94}
\begin{equation} \label{defvcsammom}
var (m_r^{\prime}) = \frac1n (\mu_{2r}^{\prime} - \mu_r^{\prime 2}), \qquad cov (\mu_r^{\prime}, \mu_q^{\prime}) = \frac1n (\mu_{r+q}^{\prime} - \mu_r^{\prime} \mu_q^{\prime})
\end{equation}
(which is an exact result) and equate them to the variance and covariance of population moments~(\ref{varmu}) and~(\ref{covmu}).  As a result we obtain the following matrix equation
\begin{equation} \label{vcmm}
\mathbf{V}_m^T = \mathbf{H} \cdot \mathbf{V}_{\theta}^T
\end{equation}
The sought $\mathbf{V}_{\theta}^T$ is defined as
\begin{equation} \label{defvtheta}
\mathbf{V}_{\theta} = (var (\theta_1), var (\theta_2), var (\theta_3), cov (\theta_1, \theta_2), cov (\theta_1, \theta_3), cov (\theta_2, \theta_3))
\end{equation}
where superscript $T$ signifies transposition.  The vector of variances and covariances of sample moments (using~(\ref{defvcsammom})) is given by:
\begin{equation}\label{vcsammom}
n \mathbf{V}_m^T = n
\left( \begin{array}{c}
var (m_1^{\prime})\\
var (m_2^{\prime})\\
var (m_3^{\prime})\\
cov (m_1^{\prime}, m_2^{\prime})\\
cov (m_1^{\prime}, m_3^{\prime})\\
cov (m_2^{\prime}, m_3^{\prime})
\end{array} \right)
= \left( \begin{array}{c}
\mu_2^{\prime} - \mu_1^{\prime 2}\\
\mu_4^{\prime} - \mu_2^{\prime 2}\\
\mu_6^{\prime} - \mu_3^{\prime 2}\\
\mu_3^{\prime} - \mu_1^{\prime} \mu_2^{\prime}\\
\mu_4^{\prime} - \mu_1^{\prime} \mu_3^{\prime}\\
\mu_5^{\prime} - \mu_2^{\prime} \mu_3^{\prime}
\end{array} \right)
\end{equation}
and can be calculated using formulas for uncorrected moments in terms of Lerch's transcendent~(\ref{compuncor}) (not shown).  The remaining component of Equation~(\ref{vcmm}) is the matrix $\mathbf{H}$, given by 
\begin{eqnarray} \label{hmatr}
\mathbf{H} & = &
\left( \begin{array}{cccc}
H_{11}^2& H_{12}^2& H_{13}^2& 2H_{11} H_{12}\\
H_{21}^2& H_{22}^2& H_{23}^2& 2H_{21} H_{22}\\
H_{31}^2& H_{32}^2& H_{33}^2& 2H_{31} H_{32}\\
H_{11} H_{21}& H_{12} H_{22}& H_{13} H_{23}& H_{11} H_{22}+H_{12} H_{21}\\
H_{11} H_{31}& H_{12} H_{32}& H_{13} H_{33}& H_{11} H_{32}+H_{12} H_{31}\\
H_{21} H_{31}& H_{22} H_{32}& H_{23} H_{33}& H_{21} H_{32}+H_{22} H_{31}
\end{array} \right. \nonumber \\
&& \left. \begin{array}{cc}
2H_{11} H_{13}& 2H_{12} H_{13}\\
2H_{21} H_{23}& 2H_{22} H_{23}\\
2H_{31} H_{33}& 2H_{32} H_{33}\\
H_{11} H_{23}+H_{13} H_{21}& H_{12} H_{23}+H_{13} H_{22}\\
H_{11} H_{33}+H_{13} H_{31}& H_{12} H_{33}+H_{13} H_{32}\\
H_{21} H_{33}+H_{23} H_{31}& H_{22} H_{33}+H_{23} H_{32}
\end{array} \right)
\end{eqnarray}
The elements of $\mathbf{H}$ 
\begin{equation} \label{helem}
H_{ij} = \frac{\partial \mu_i^{\prime}}{\partial \theta_j}
\end{equation}
can be expressed in terms of Lerch's transcendent (Appendix \ref{AppA}).

Finally, the vector of variances and covariances of estimates $\bar{\theta}=(\bar{z},\bar{s},\bar{v})$ is obtained as the solution of Equation~(\ref{vcmm}) by using~(\ref{vcsammom}), (\ref{hmatr}) and~(\ref{helem1}):
\begin{equation} \label{mmvar}
\mathbf{V}_{\theta}^T = \mathbf{H}^{-1} \cdot \mathbf{V}_{m}^T
\end{equation}

\subsection{Method of maximum likelihood}

Maximum likelihood (ML) estimators for the Lerch distribution can be obtained in a relatively straightforward way.  The likelihood function is
\begin{equation} \label{like}
L(\theta) = \prod_{i=1}^{n} Pr(X=x_i) = c^n \frac{z^{\sum_{i=1}^{n} x_i}}{\left( \prod_{i=1}^{n} (v+x_i) \right)^s}
\end{equation}
It is more convenient to work with the log-likelihood function which is obtained by taking the logarithm of~(\ref{like})
\begin{equation} \label{loglike}
\log L(\theta) = -n \log \Phi(z,s,v) + \log z \sum_{i=1}^{n} x_i - s \sum_{i=1}^{n} \log (v+x_i)
\end{equation}
ML estimates are obtained by maximizing the likelihood function $L(\theta)$~(\ref{like}) or equivalently minimizing the log-likelihood $\log L(\theta)$~(\ref{loglike}), which is done by solving the following system of equations for derivatives
\begin{eqnarray} \label{loglikeeq}
\frac{\partial \log L}{\partial z} & = & -n \frac{\partial \log \Phi}{\partial z} + \frac{\sum_{i=1}^{n} x_i}{z} = 0 \nonumber \\
\frac{\partial \log L}{\partial s} & = & -n \frac{\partial \log \Phi}{\partial s} - \sum_{i=1}^{n} \log (v+x_i) = 0 \nonumber \\
\frac{\partial \log L}{\partial v} & = & -n \frac{\partial \log \Phi}{\partial v} - s \sum_{i=1}^{n} \frac1{v+x_i} = 0
\end{eqnarray}
Substituting the partial derivatives of Lerch's transcendent from Appendix \ref{AppA} into Equations~(\ref{loglikeeq}) we obtain the following system of equations
\begin{eqnarray} \label{mle}
n^{-1} \sum_{i=1}^{n} x_i & = & \sum_j f_j x_j = \hat{z} \frac{\Phi(\hat{z},\hat{s}-1,\hat{v}+1) - \hat{v} \Phi(\hat{z},\hat{s},\hat{v}+1)}{\Phi(\hat{z},\hat{s},\hat{v})} \nonumber \\
n^{-1} \sum_{i=1}^{n} \log (\hat{v}+x_i) & = & \sum_j f_j \log(\hat{v}+x_j) = \frac{\sum_{x=0}^{\infty} \frac{\log (\hat{v}+x) \hat{z}^x}{(\hat{v}+x)^{\hat{s}}}}{\Phi(\hat{z},\hat{s},\hat{v})} \nonumber \\
n^{-1} \sum_{i=0}^{n} \frac1{\hat{v}+x_i} & = & \sum_j \frac{f_j}{\hat{v}+x_j} = \frac{\Phi(\hat{z},\hat{s}+1,\hat{v})}{\Phi(\hat{z},\hat{s},\hat{v})}
\end{eqnarray}

To obtain the ML estimates $\hat{z}$, $\hat{s}$ and $\hat{v}$ we have to solve numerically Equations~(\ref{mle}).

Asymptotic variances and covariances of the ML estimates~(\ref{defvtheta}) are obtained as the inverse of Fisher's information matrix whose $ij$th element is the expectation of the second partial derivative of the log-likelihood function with respect to the parameters, evaluated at the ML estimates: 
\begin{equation} \label{ielem}
n I_{ij} = E\left(- \frac{\partial^2 \log L(\hat{\theta})}{\partial \theta_i \partial \theta_j} \right)
\end{equation}
The information matrix is then defined as
\begin{equation} \label{fim}
\mathbf{I}= n 
\left( \begin{array}{ccc}
I_{11} & I_{12} & I_{13}\\
I_{12} & I_{22} & I_{23}\\
I_{13} & I_{23} & I_{33}
\end{array} \right)
\end{equation}
The elements~(\ref{ielem}) of the information matrix can be obtained by straightforward differentiation and expressed in terms of Lerch's transcendent (Appendix \ref{AppA}).

Finally, the elements of the vector of variances and covariances $\mathbf{V}_{\theta}$~(\ref{defvtheta}) can be arranged as a variance-covariance matrix
\begin{equation} \label{vcmatr}
\mathbf{V} = 
\left( \begin{array}{ccc}
var (z)& cov(z,s)& cov(z,v)\\
cov(z,s)& var (s)& cov(s,v)\\
cov(z,v)& cov(s,v)& var (v)
\end{array} \right)
\end{equation}
which is calculated using Equations~(\ref{fim}) and~(\ref{ielem1}):
\begin{equation} \label{mlevar}
\mathbf{V} = \mathbf{I}^{-1}
\end{equation}

With both fitting methods, goodness-of-fit can be determined by using the classical Pearson Chi-squared statistic $X^2 = n \sum_{j=1}^M (f_j - p_j)^2/ p_j$, where $M$ is the number of observed classes, $f_j$ is the observed frequency for the $j$th class, and $p_j(\theta)$ is the fitted frequency (p.m.f.) for the $j$th class \citep{AgSt86}.  Even though we are really testing the composite hypothesis (we estimate parameters from data), asymptotically $X^2$ is distributed as $\chi^2(M-p-1)$, where $p=3$ is the number of parameters estimated.  Minimization of the Chi-squared statistic can also be used to estimate parameters of the Lerch distribution.

\section{Examples}

In this section we consider several examples of fitting the Lerch distribution to over- and under-dispersed data that arise in a variety of counting processes, and to rank-abundance ecological data.  These data sets have been modeled using Poisson, Poisson mixtures, generalized Poisson, adjusted generalized Poisson, and Zipf-Mandelbrot distributions.  Our results show that the Lerch distribution provides better goodness-of-fit and in two cases a successful fit when the alternatives fail to fit in the whole range of data.

\begin{table} 
\begin{center}
\begin{tabular}{llll}
\hline
Counts & Observed & Gen. Poisson & Lerch \\
\hline
0 & 28 & 26.1127 & 29.2839 \\
1 & 28 & 23.6448 & 21.153 \\ 
2 & 14 & 18.095 & 15.6055 \\ 
3 & 11 & 13.3871 & 11.7173 \\ 
4 & 8 & 9.86865 & 8.93021 \\ 
5 & 11 & 7.31252 & 6.89379 \\ 
6 & 2 & 5.46003 & 5.38124 \\
7 & 3 & 4.10969 & 4.24166 \\ 
8 & 3 & 3.11713 & 3.37228 \\ 
9 & 3 & 2.38107 & 2.70168 \\ 
10 & 3 & 1.83053 & 2.1791 \\ 
11 & 2 & 1.41548 & 1.76882 \\
12 & 0 & 1.10029 & 1.44369 \\
13 & 1 & 0.859339 & 1.18433 \\ 
14 & 2 & 0.67404 & 0.976077 \\ 
15 & 1 & 0.530761 & 0.807877 \\ 
16 & 0 & 0.419425 & 0.671287 \\ 
17 & 2 & 0.332522 & 0.559812 \\
\hline
\end{tabular}
\end{center}
\caption{Distribution of sowbugs under boards.  The observations are frequencies of sowbugs per board, frequency classes are numbers of sowbugs per board.  Observed frequencies are taken from \citep{Co46}. Predicted frequences come from the best-fit generalized Poisson and Lerch distributions.  See text for details.} \label{tab:sowbugs}
\end{table}

The first data set that we analyze is the overdispersed data on the numbers of sowbugs \emph{Trachelipus rathkei} found under boards (columns 1 and 2 in Table~\ref{tab:sowbugs}).  The data were obtained by \cite{Co46} in studies of the distribution of different cryptozoa species within areas of their natural habitat.  Cole found that the distribution of spiders under wooden boards scattered in the area follows a Poisson distribution, which indicates a random distribution of individuals and their unsocial behavior.  In contrast, the distribution of sowbugs exhibited properties of the contagious distribution and could not be fit with a Poisson distribution. \cite{JaKeSc79} successfully fit the sowbug data set with the generalized Poisson model that has the p.f. \citep{Co89}
\begin{equation} \label{genp}
p_x = \left\{ 
\begin{array}{ll}
\frac{\theta (\theta + x \lambda)^{x-1} \mathrm{e}^{-\theta -x \lambda}}{x!}, & x=0,1,2,\dots \\
0,& x>\left[-\frac{\theta}{\lambda}\right] \quad \mathrm{when} \quad \lambda<0
\end{array} \right.
\end{equation}
where $[\cdot]$ signifies the integer part operation, with parameters $\theta=1.5416$ and $\lambda=0.5321$ (see column 3 in Table~\ref{tab:sowbugs} for calculated frequencies). We grouped classes 6 and 7, 8 and 9, 10 and 11, and 12 through 17.  With this grouping we obtained the Pearson statistic $X^2=9.3089$\footnote{All numbers in this section have been rounded to five significant digits}, with a p-value of 0.231232 for 7 degrees of freedom (d.f.).

Here we fit this data set with the Lerch distribution.  To fit the Lerch distribution, we use the same grouping, and minimize Pearson $X^2$ statistic.  The best-fit parameters are $z=0.913315$, $s=2.37621$, $v=9.63785$, and the achieved minimum $X^2=7.69169$ (see column 4 in Table~\ref{tab:sowbugs} for calculated frequencies).  The p-value is 0.261572 for 6 d.f.  We see that both generalized Poisson and Lerch fits are acceptable and goodness-of-fit is better with the Lerch distribution.

\begin{table} 
\begin{center}
\begin{tabular}{llll}
\hline
Counts & Observed & Adj. Gen. Poisson & Lerch \\
\hline
0 & 162 & 162.004 & 161.906 \\
1 & 267 & 264.773 & 266.73 \\
2 & 271 & 268.751 & 264.789 \\
3 & 185 & 194.406 & 192.091 \\
4 & 111 & 112.384 & 112.56 \\
5 & 61 & 55.2168 & 56.4979 \\
6 & 27 & 23.9531 & 25.217 \\
7 & 8 & 9.4128 & 10.2649 \\
8 & 3 & 3.41263 & 3.87964 \\
9 & 1 & 1.15712 & 1.37948 \\
\hline
\end{tabular}
\end{center}
\caption{Distribution of death notices in the London ``Times''.  The observations are frequencies of death notices published per day, frequency classes are numbers of notices published per day.  Observed frequencies are taken from \citep{Ha69}. Predicted frequences come from the best-fit adjusted generalized Poisson and Lerch distributions.  See text for details.} \label{tab:death}
\end{table}

Another way to model data that is over- or under-dispersed is to use an adjusted distribution that would correct for larger or smaller number of zeros \citep{GuGuTr96}.  Consider data on the numbers of death notices for women 80 years and older that had been published in the London ``Times'' every day for three consecutive years (columns 1 and 2 in Table~\ref{tab:death}) \citep{Ha69}.  The data were examined for differences of death rates between winter and summer months and thus were originally fitted by a mixture of two Poisson distributions.  \cite{GuGuTr96} observed that the data set might be zero deflated and fit an adjusted generalized Poisson model in which they introduced an additional parameter $\phi$ to describe fewer zeros than predicted by the generalized Poisson:
\begin{equation} \label{genpz}
p_x  = \left\{ \begin{array}{ll}
\phi + (1-\phi) \mathrm{e}^{-\theta}, & x=0 \\ 
\frac{(1 + x \alpha)^{x-1}}{x!} \frac{(\theta \mathrm{e}^{-\alpha \theta})^x}{\mathrm{e}^\theta}, & x=1,2,\dots \\
\end{array} \right.
\end{equation}
P.m.f.~(\ref{genp}) is given here in a so-called restricted form \citep{Co89} as in \citep{GuGuTr96}.  The best-fit parameters of the adjusted generalized Poisson model are $\theta=2.038$\footnote{Here we corrected an apparent misprint in Table~3 and Figure~4 of \citep{GuGuTr96} that read $\theta=1.2038$}, $\alpha=0.03639$ and $\phi=0.02015$ (see column 3 in Table~\ref{tab:death} for calculated frequencies).  The Pearson statistic is $X^2=1.9379$, with a p-value of 0.7472 for 4 d.f.

Since the Lerch distribution accounts for the zero class naturally, we fitted it to the death notices data by minimizing the Pearson $X^2$ statistic.  The best-fit parameters are $z=0.189628$, $s=-7.10717$, $v=2.81275$ (see column 4 in Table~\ref{tab:death} for calculated frequencies).  Minimum Pearson statistic is $X^2=1.23938$, with a p-value of 0.871573 for 4 d.f.   Note that since the last three classes were grouped in \citep{Ha69} and \citep{GuGuTr96}, we grouped them here as well to fit the Lerch distribution and calculate the minimum $X^2$.  Goodness-of-fit is better with the Lerch distribution than with the adjusted generalized Poisson model.

\begin{table} 
\begin{center}
\begin{tabular}{lllll}
\hline
Counts & 0 & 1 & 2 & 3  \\
Observed & 5 & 68 & 88 & 32 \\
Gen. Poisson & 8.67105 & 57.5966 & 97.4296 & 29.5181 \\
Lerch & 2.79364 & 67.0555 & 93.2649 & 26.8808 \\
\hline
\end{tabular}
\end{center}
\caption{Distribution of bean weevil eggs on beans.  The observations are frequencies of eggs oviposited by bean weevil per bean, frequency classes are numbers of eggs found on a bean.  Observed frequencies are taken from \citep{Mi75}. Predicted frequences come from the best-fit generalized Poisson and Lerch distributions.  See text for details.} \label{tab:beans}
\end{table}

Now we turn to underdispersed data on the numbers of eggs oviposited by bean weevil \emph{Callosobrachus maculatus} on beans (rows 1 and 2 in Table~\ref{tab:beans}) \citep{Mi75}.   The data were obtained in studies of oviposition tactics by bean weevils with the mean number of eggs per bean 1.8.  Again, a simple Poisson distribution did not fit well, suggesting that egged and unegged beans are not equally attractive for weevils seeking to oviposit an egg.  \cite{JaKeSc79} fit the generalized Poisson distribution with~(\ref{genp}) to the weevil data with parameters $\theta=3.1027$ and $\lambda=-0.7612$ (see row 3 in Table~\ref{tab:beans} for calculated frequencies).  For comparison with the Lerch distribution, for which we would have zero d.f. of the Chi-squared distribution, we use the sum of squared deviations between the p.f. and empirical frequencies as a goodness-of-fit measure.  For the generalized Poisson distribution the sum is 0.00581991.

We fit the Lerch distribution to the weevil data by minimizing the Pearson $X^2$ statistic.  The best-fit parameters are $z=0.00116201$, $s=-24.9577$, $v=2.04499$ (see row 4 in Table~\ref{tab:beans} for calculated frequencies).  The sum of squared deviations is 0.00160233, and the goodness-of-fit is again better with the Lerch distribution.

Finally, we consider rank-abundance data on biotic compartments in Lake Yunoko, Japan  (rows 1 and 2 in Table~\ref{tab:yunoko}) \citep{Ao95}.  The Zipf-Mandelbrot distribution as a model for frequencies of ranked species was suggested for use in ecology by \cite{Fr85}. \cite{Ao95} observed that neither Zipf nor Zipf-Mandelbrot models can fit the whole range of data and suggested disregarding the most abundant species (fish) from the data, which could then be fit with the Zipf distribution.  \cite{Ao95} gives the anthropogenic source of fish in the lake as a reason for disregarding the one data point.  However, there are reasons to believe that different biotic compartments influence each other and, in particular the disregarded fish, might influence abundances of phytoplankton and other compartments. Thus, it is of interest to fit the whole range of data.  The Lerch distribution provides a suitable model, being a generalization of the Zipf-Mandelbrot distribution.  There are only six classes in the data starting with one.  The observed abundancies are calculated from the data on standing crop in compartments, expressed in grams of carbon per m$^2$ \citep{Ao95}, by dividing these values by the total crop.  We fit the doubly truncated Lerch distribution by minimizing the Pearson $X^2$ statistic.  The best-fit parameters are $z=0.219158$, $s=-0.214704$, $v=-0.998437$ (see row 3 in Table~\ref{tab:yunoko} for calculated frequencies); the Pearson statistic is $X^2=0.0259897$, with a p-value of 0.987089 for 2 d.f.  The fit is obviously very good.

\begin{table} 
\begin{center}
\begin{tabular}{lllllll}
\hline
Rank & 1 & 2 & 3 & 4 & 5 & 6 \\
Observed & 0.46798 & 0.428571 & 0.0738916 & 0.0152709 & 0.00837438 & 0.00591133 \\
Lerch & 0.460902 & 0.404575 & 0.102876 & 0.0245955 & 0.00573359 & 0.00131821 \\
\hline
\end{tabular}
\end{center}
\caption{Distribution of ranked biotic compartments in Lake Yunoko.  The observations are relative abundances of the biotic compartments, classes are ranks of compartments (rank one is the most abundant, etc.) Observed frequencies are taken from \citep{Ao95}.  Predicted frequences come from the best-fit Lerch distribution.  See text for details.} \label{tab:yunoko}
\end{table}

In summary, the above examples demonstrate that the Lerch distribution can provide a useful model with better goodness-of-fit for over- and under-dispersed data, for which modifications of the Poisson distribution or the reduced Lerch (Good) model perform less well or even fail, as well as for some rank-abundance ecological data, for which the reduced Lerch (Zipf-Mandelbrot) model does not fit.  The ultimate reason is that the Lerch distribution is a generalization of the Good and Zipf-Mandelbrot distributions (which were originally proposed as models for these problems) and is endowed with greater flexibility to accomodate the odd distributional shapes often observed in real problems.  A comment on the method of parameter estimation is in order.  We found in practice that the most convenient and accurate method of parameter estimation is minimization of the Pearson $X^2$ statistic, which works well for data even with a few frequency classes.  The moment and maximum likelihood methods are more difficult to use in such situations because of higher variance of sample moments and expectations in a small sample setting.  An added difficulty is the slow convergence of MM and ML equations that involve Lerch's transcendent.  These difficulties are not specific to the Lerch distribution and are in fact common for other distributions that are based on special functions.  However, we found that these methods perform well with more rich data sets (not shown here).

\appendix

\section{Appendix: variance-covariance matrices for method of moments and method of maximum likelihood estimators}
\label{AppA}

To calculate the elements of matrices $\mathbf{H}$ and $\mathbf{I}$, we need the derivatives of Lerch's transcendent, which can be calculated by straightforward differentiation term by term:
\begin{eqnarray} \label{lerchderiv}
\frac{\partial \Phi(z,s,v)}{\partial z} & = & \Phi(z,s-1,v+1) - v \Phi(z,s,v+1) \nonumber \\ 
\frac{\partial \Phi(z,s,v)}{\partial s} & = & -\sum_{n=0}^{\infty} \frac{\log (v+n) z^n}{(v+n)^s} \nonumber \\
\frac{\partial \Phi(z,s,v)}{\partial v} & = & -s \Phi(z,s+1,v)
\end{eqnarray}
Note that the infinite sum in the second equation~(\ref{lerchderiv}) converges, because $\log (v+x)$ grows slower than $(v+x)$, which converges to $\Phi(z,s-1,v)$.

Then the quantities $H_{ij}$~(\ref{helem}), which are the elements of  $\mathbf{H}$~(\ref{hmatr}) are calculated as follows
\begin{eqnarray} \label{helem1}
H_{11} & = & \Phi^{-2}(z,s,v) ( \Phi(z,s-2,v+1) \Phi(z,s,v) - \nonumber \\
&& \Phi(z,s-1,v+1) \Phi(z,s-1,v) + v (\Phi(z,s,v+1) \Phi(z,s-1,v) - \nonumber \\
&& \Phi(z,s-1,v+1) \Phi(z,s,v))) \nonumber \\
H_{12} & = & \Phi^{-2}(z,s,v) \left( \Phi(z,s-1,v) \sum_{n=0}^{\infty} \frac{\log (v+n) z^n}{(v+n)^s} - \right. \nonumber \\
&& \left. \Phi(z,s,v) \sum_{n=0}^{\infty} \frac{\log (v+n) z^n}{(v+n)^{s-1}}  \right) \nonumber \\
H_{13} & = & \Phi^{-2}(z,s,v) ( (1-s) \Phi^2(z,s,v) + s \Phi(z,s-1,v) \Phi(z,s+1,v) ) - 1\nonumber \\
H_{21} & = & -2v H_{11} + \Phi^{-2}(z,s,v) ( \Phi(z,s-3,v+1) \Phi(z,s,v) - \nonumber \\
&& \Phi(z,s-1,v+1) \Phi(z,s-2,v) + v (\Phi(z,s,v+1) \Phi(z,s-2,v) - \nonumber \\
&& \Phi(z,s-2,v+1) \Phi(z,s,v))) \nonumber \\
H_{22} & = & -2v H_{12} + \Phi^{-2}(z,s,v) \left( \Phi(z,s-2,v) \sum_{n=0}^{\infty} \frac{\log (v+n) z^n}{(v+n)^s} - \right. \nonumber \\
&& \left. \Phi(z,s,v) \sum_{n=0}^{\infty} \frac{\log (v+n) z^n}{(v+n)^{s-2}} \right) \nonumber \\
H_{23} & = & -2v H_{13} - 2 \frac{\Phi(z,s-1,v)}{\Phi(z,s,v)} + \nonumber \\
&& \Phi^{-2}(z,s,v) ( (2-s) \Phi(z,s-1,v) \Phi(z,s,v) + \nonumber \\
&& s \Phi(z,s+1,v) \Phi(z,s-2,v) ) \nonumber \\
H_{31} & = & -3v^2 H_{11} -3v H_{21} + \Phi^{-2}(z,s,v) ( \Phi(z,s-4,v+1) \Phi(z,s,v) - \nonumber \\
&& \Phi(z,s-1,v+1) \Phi(z,s-3,v) + v (\Phi(z,s,v+1) \Phi(z,s-3,v) - \nonumber \\
&& \Phi(z,s-3,v+1) \Phi(z,s,v))) \nonumber \\
H_{32} & = & -3v^2 H_{12} -3v H_{22} + \Phi^{-2}(z,s,v) \left( \Phi(z,s-3,v) \sum_{n=0}^{\infty} \frac{\log (v+n) z^n}{(v+n)^s} - \right. \nonumber \\
&& \left. \Phi(z,s,v) \sum_{n=0}^{\infty} \frac{\log (v+n) z^n}{(v+n)^{s-3}} \right) \nonumber \\
H_{33} & = & -3v^2 H_{13} -3v H_{23} -3 \frac{\Phi(z,s-2,v)}{\Phi(z,s,v)} + \nonumber \\
&& \Phi^{-2}(z,s,v) ( (3-s) \Phi(z,s-2,v) \Phi(z,s,v) + \nonumber \\
&& s \Phi(z,s+1,v) \Phi(z,s-3,v) )
\end{eqnarray}

The elements $I_{ij}$~(\ref{ielem}) of the information matrix~(\ref{fim}) are calculated as follows
\begin{eqnarray} \label{ielem1}
I_{11} & = & \Phi^{-2}(z,s,v) ( \Phi(z,s,v) ( \Phi(z,s-2,v+2) - \nonumber \\
&& (2 v+1) \Phi(z,s-1,v+2) + 
v (v+1) \Phi(s,z,v+2)) - \nonumber \\
&& (\Phi(z,s-1,v) - v \Phi(z,s,v+1))^2) + 
\frac1{z^2} \left( \frac{\Phi(z,s-1,v)}{\Phi(z,s,v)}-v \right) \nonumber \\
I_{12} & = & \Phi^{-2}(z,s,v) \left( (\Phi(z,s-1,v+1) - v \Phi(z,s,v+1)) \times \right. \nonumber \\
&& \left. \sum_{n=0}^{\infty} \frac{\log (v+n) z^n}{(v+n)^s} - \Phi(z,s,v) \sum_{n=0}^{\infty} \frac{\log (v+n+1) (n+1) z^n}{(v+n+1)^s} \right) \nonumber \\
I_{13} & = & s \Phi^{-2}(z,s,v) (\Phi(z,s,v) (v \Phi(z,s+1,v+1) - \Phi(z,s,v+1)) + \nonumber \\
&& \Phi(z,s+1,v) (\Phi(z,s-1,v+1) - v \Phi(z,s,v+1))) \nonumber \\
I_{22} & = & \Phi^{-2}(z,s,v) \left( \Phi(z,s,v) \sum_{n=0}^{\infty} \frac{\log^2 (v+n) z^n}{(v+n)^s} - \right. \nonumber \\
&& \left. \left( \sum_{n=0}^{\infty} \frac{\log (v+n) z^n}{(v+n)^s} \right)^2 \right) \nonumber \\
I_{23} & = & \Phi^{-2}(z,s,v) \left( \Phi(z,s,v) \sum_{n=0}^{\infty} \frac{(s \log(v+n)-1) z^n}{(v+n)^{s+1}} - \right. \nonumber \\
&& \left. s \Phi(z,s+1,v) \sum_{n=0}^{\infty} \frac{\log (v+n) z^n}{(v+n)^s} \right) + \frac{\Phi(z,s+1,v)}{\Phi(z,s,v)} \nonumber \\
I_{33} & = & s^2 \Phi^{-2}(z,s,v) (\Phi(z,s,v) \Phi(z,s+2,v) - \Phi^2(z,s+1,v))
\end{eqnarray}
When calculating expectations~(\ref{ielem}), we use the following asymptotic results:
\begin{eqnarray} \label{asymexp}
E\left(\frac1n \sum_{i=1}^n x_i\right) & = & \mu_1^{\prime} \nonumber \\
E\left(\frac1n \sum_{i=1}^n \frac1{v+x_i}\right) & = & \frac{\Phi(z,s+1,v)}{\Phi(z,s,v)} \nonumber \\
E\left(\frac1n \sum_{i=1}^n \frac1{(v+x_i)^2}\right) & = & \frac{\Phi(z,s+2,v)}{\Phi(z,s,v)}
\end{eqnarray}

\bibliography{REFS}